\documentclass[12pt]{article}
\usepackage{amssymb,amsmath}
\usepackage{hyperref}
\usepackage{graphicx}
\usepackage{color}

\begin{document}

\title{\LARGE\bf A simple identity for derivatives of the arctangent function}

\author{
\normalsize\bf S. M. Abrarov\footnote{\scriptsize{Dept. Earth and Space Science and Engineering, York University, Toronto, Canada, M3J 1P3.}}\, and B. M. Quine$^{*}$\footnote{\scriptsize{Dept. Physics and Astronomy, York University, Toronto, Canada, M3J 1P3.}}}

\date{May 10, 2016}
\maketitle

\begin{abstract}
We present an identity for the derivatives of the arctangent function as an alternative to the Adegoke--Layeni--Lampret formula. We show that algorithmic implementation of the proposed identity can significantly accelerate the computation since this approach requires no symbolic programming in determination of the derivatives for the arctangent function.
\vspace{0.25cm}
\\
\noindent {\bf Keywords:} arctangent function, Adegoke--Layeni--Lampret formula, numerical integration
\vspace{0.25cm}
\end{abstract}

\section{Introduction}

The arctangent function that can be represented in form of the integral
\begin{equation}\label{eq_1}
\arctan \left( x \right) = \int\limits_0^1 {\frac{x}{{1 + {x^2}{t^2}}}dt}
\end{equation}
is widely used in many applications. Historically, many identities for the arctangent function were discovered centuries ago \cite{Lehmer1938, Beckmann1976, Abeles1993, Berggren2004, Goyanes2010}. Despite this, finding new identities of the arctangent function is still topical according to the modern literature \cite{Calcut1999, Adegoke2010, Chen2010, Lampret2011, Sofo2012}.
Recently Lampret \cite{Lampret2011} reported a remarkable formula for the arctangent function 
\begin{equation}\label{eq_2}
\frac{{{d^m}}}{{d{t^m}}}\arctan \left( t \right) = {\rm{sg}}{{\rm{n}}^{m - 1}}\left( { - t} \right)\frac{{\left( {m - 1} \right)!}}{{{{\left( {1 + {t^2}} \right)}^{m/2}}}}\sin \left( {m\arcsin \left( {\frac{1}{{\sqrt {1 + {t^2}} }}} \right)} \right),
\end{equation}
where
$$
{\rm{sgn}}\left( t \right) = \left\{
\begin{aligned}
1, & \qquad t \ge 0\\
 - 1, & \qquad t < 0
\end{aligned}
\right.
$$
is the signum function. Notably, Lampret in his work referred to variation of this formula without inclusion of factor $\rm{sgn}^{m-1} \left( t \right)$ that had been given earlier by Adegoke and Layeni \cite{Adegoke2010}. Therefore, this equation may be regarded as the Adegoke--Layeni--Lampret formula.

The Adegoke--Layeni--Lampret formula \eqref{eq_2} is indeed elegant. However, application of this formula causes challenges in algorithmic implementation. Although in the Adegoke--Layeni--Lampret formula \eqref{eq_2} the multiplier
\[
\sin \left( {m\arcsin \left( {\frac{1}{{\sqrt {1 + {t^2}} }}} \right)} \right)
\]
changes in a periodical sequence with increasing order $m$ in derivative of the arctangent function, this multiplier is not an elementary function. As a result, it cannot be used for general application in numerical integration. In this work we derive an identity for derivatives of the arctangent function in form of a rational function and, as an example, we show that due to its simple representation the proposed identity can be used for efficient computation of the arctangent function and the constant pi.

\section{Derivation by induction procedure}

Consider the following function
\begin{equation}\label{eq_3}
\frac{1}{{1 - {u^2}}}.
\end{equation}
It is not difficult to see that this function can be represented in form
$$
\frac{1}{{1 - {u^2}}} = {\left( { - 1} \right)^0}\frac{{0!}}{2}\left( {\frac{1}{{{{\left( {u + 1} \right)}^1}}} - \frac{1}{{{{\left( {u - 1} \right)}^1}}}} \right).
$$
Taking derivative of this function results in
$$
\begin{aligned}
\frac{d}{{du}}\left( {\frac{1}{{1 - {u^2}}}} \right) &= {\left( { - 1} \right)^0}\frac{{0!}}{2}\frac{d}{{du}}\left( {\frac{1}{{{{\left( {u + 1} \right)}^1}}} - \frac{1}{{{{\left( {u - 1} \right)}^1}}}} \right) \\
&= {\left( { - 1} \right)^1}\frac{{1!}}{2}\left( {\frac{1}{{{{\left( {u + 1} \right)}^2}}} - \frac{1}{{{{\left( {u - 1} \right)}^2}}}} \right).
\end{aligned}
$$
Taking the second and the third derivatives of the function \eqref{eq_3} leads to
$$
\begin{aligned}
\frac{{{d^2}}}{{d{u^2}}}\left( {\frac{1}{{1 - {u^2}}}} \right) &= {\left( { - 1} \right)^1}\frac{{1!}}{2}\frac{d}{{du}}\left( {\frac{1}{{{{\left( {u + 1} \right)}^2}}} - \frac{1}{{{{\left( {u - 1} \right)}^2}}}} \right) \\
&= {\left( { - 1} \right)^2}\frac{{2!}}{2}\left( {\frac{1}{{{{\left( {u + 1} \right)}^3}}} - \frac{1}{{{{\left( {u - 1} \right)}^3}}}} \right)
\end{aligned}
$$
and
$$
\begin{aligned}
\frac{{{d^3}}}{{d{u^3}}}\left( {\frac{1}{{1 - {u^2}}}} \right) &= {\left( { - 1} \right)^2}\frac{{2!}}{2}\frac{d}{{du}}\left( {\frac{1}{{{{\left( {u + 1} \right)}^3}}} - \frac{1}{{{{\left( {u - 1} \right)}^3}}}} \right) \\
&= {\left( { - 1} \right)^3}\frac{{3!}}{2}\left( {\frac{1}{{{{\left( {u + 1} \right)}^4}}} - \frac{1}{{{{\left( {u - 1} \right)}^4}}}} \right),
\end{aligned}
$$
respectively. Repeating step-by-step the same procedure to higher orders of derivatives we can find such a way by induction that
$$
\begin{aligned}
\frac{{{d^m}}}{{d{u^m}}}\left( {\frac{1}{{1 - {u^2}}}} \right) &= \frac{d}{{du}}\left( {\frac{{{d^{m - 1}}}}{{d{u^{m - 1}}}}\left( {\frac{1}{{1 - {u^2}}}} \right)} \right) \\
&= {\left( { - 1} \right)^{m - 1}}\frac{{\left( {m - 1} \right)!}}{2}\frac{d}{{du}}\left( {\frac{1}{{{{\left( {u + 1} \right)}^m}}} - \frac{1}{{{{\left( {u - 1} \right)}^m}}}} \right)
\end{aligned}
$$
since in each consecutive step we {\it{de facto}} always apply the same rule
$$
\frac{d}{{du}}\left( {\frac{1}{{{{\left( {u \pm 1} \right)}^m}}}} \right) =  - \frac{m}{{{{\left( {u \pm 1} \right)}^{m + 1}}}}, \quad m \ge 1.
$$
Consequently, taking this into consideration and observing that for the zeroth derivative
$$
\frac{{{d^0}}}{{d{u^0}}}\left( {\frac{1}{{1 - {u^2}}}} \right) = \frac{1}{{1 - {u^2}}} = {\left( { - 1} \right)^0}\frac{{0!}}{2}\left( {\frac{1}{{{{\left( {u + 1} \right)}^1}}} - \frac{1}{{{{\left( {u - 1} \right)}^1}}}} \right)
$$
we can conclude that any order of the derivatives can be found according to the equation
\begin{equation}\label{eq_4}
\frac{{{d^m}}}{{d{u^m}}}\left( {\frac{1}{{1 - {u^2}}}} \right) = {\left( { - 1} \right)^m}\frac{{m!}}{2}\left( {\frac{1}{{{{\left( {u + 1} \right)}^{m + 1}}}} - \frac{1}{{{{\left( {u - 1} \right)}^{m + 1}}}}} \right), \quad m \ge 0.
\end{equation}

Since
$$
\frac{{{d^m}}}{{d{t^m}}}\left( {\frac{1}{{1 + {t^2}}}} \right) = \frac{{{d^m}}}{{d{t^m}}}\left( {\frac{1}{{1 - {{\left( {it} \right)}^2}}}} \right) = {i^m}\frac{{{d^m}}}{{d{{\left( {it} \right)}^m}}}\left( {\frac{1}{{1 - {{\left( {it} \right)}^2}}}} \right)
$$
at $u = i t$ from the equation \eqref{eq_4} we immediately obtain a very useful identity
\begin{equation}\label{eq_5}
\frac{{{d^m}}}{{d{t^m}}}\left( {\frac{1}{{1 + {t^2}}}} \right) = {\left( { - i} \right)^m}\frac{{m!}}{2}\left( {\frac{1}{{{{\left( {it + 1} \right)}^{m + 1}}}} - \frac{1}{{{{\left( {it - 1} \right)}^{m + 1}}}}} \right), \qquad m \ge 0.
\end{equation}
Due to relation
$$
\int {\frac{1}{{1 + {t^2}}}dt = } \arctan \left( t \right)
$$
we can write
\begin{equation}\label{eq_6}
\frac{{{d^m}}}{{d{t^m}}}\arctan \left( t \right) = \frac{{{d^{m - 1}}}}{{d{t^{m - 1}}}}\left( {\frac{1}{{1 + {t^2}}}} \right).
\end{equation}
Consequently, from the equation \eqref{eq_5} it follows that the derivatives for the arctangent function can be represented in form
$$
\frac{{{d^m}}}{{d{t^m}}}\arctan \left( t \right) = {\left( { - i} \right)^{m - 1}}\frac{{\left( {m - 1} \right)!}}{2}\left( {\frac{1}{{{{\left( {it + 1} \right)}^m}}} - \frac{1}{{{{\left( {it - 1} \right)}^m}}}} \right), \qquad m \ge 1
$$
or
\begin{equation}\label{eq_7}
\frac{{{d^m}}}{{d{t^m}}}\arctan \left( t \right) = \frac{{{{\left( { - 1} \right)}^m}\left( {m - 1} \right)!}}{{2i}}\left( {\frac{1}{{{{\left( {t + i} \right)}^m}}} - \frac{1}{{{{\left( {t - i} \right)}^m}}}} \right), \qquad m \ge 1.
\end{equation}
We can see now from the identity \eqref{eq_7} that the derivatives of the arctangent function can also be expressed as a simple rational function as an alternative to the Adegoke--Layeni--Lampret formula \eqref{eq_2} discussed above.

Making change of the variable in the equation \eqref{eq_7} as $t \to xt$ leads to
$$
\frac{{{\partial ^m}}}{{\partial {{\left( {xt} \right)}^m}}}\arctan \left( {xt} \right) = \frac{{{{\left( { - 1} \right)}^m}\left( {m - 1} \right)!}}{{2i}}\left( {\frac{1}{{{{\left( {xt + i} \right)}^m}}} - \frac{1}{{{{\left( {xt - i} \right)}^m}}}} \right), \quad m \ge 1
$$
or
\small
\begin{equation}\label{eq_8}
\frac{{{\partial ^m}}}{{\partial {t^m}}}\arctan \left( {xt} \right) = \frac{{{{\left( { - 1} \right)}^m}\left( {m - 1} \right)!{x^m}}}{{2i}}\left( {\frac{1}{{{{\left( {xt + i} \right)}^m}}} - \frac{1}{{{{\left( {xt - i} \right)}^m}}}} \right), \quad m \ge 1.
\end{equation}
\normalsize

\section{Applications to the arctangent function}

In our recent publication we derived the equation that can be applied for highly accurate numerical integration \cite{Abrarov2016}
\begin{equation}\label{eq_9}
\int\limits_0^1 {f\left( t \right)dt}  = \mathop {\lim }\limits_{L \to \infty } \sum\limits_{\ell  = 1}^L {\sum\limits_{m = 0}^M {\frac{{{{\left( { - 1} \right)}^m} + 1}}{{{{\left( {2L} \right)}^{m + 1}}\left( {m + 1} \right)!}}} } {\left. {{f^{\left( m \right)}}\left( t \right)} \right|_{t = \frac{{\ell  - 1/2}}{L}}},
\end{equation}
where $f\left( t \right)$ is a differentiable function in the interval $t \in \left[ {0,1} \right]$. This equation can be rearranged in a more simplified form. In particular, since
$$
{\left( { - 1} \right)^m} + 1 = \left\{ \begin{aligned}
2, &\qquad m = \left\{ {0,2,4, \ldots } \right\}\\
0, &\qquad m = \left\{ {1,3,5, \ldots } \right\}
\end{aligned} \right.
$$
we can recast the limit \eqref{eq_9} as
\begin{equation}\label{eq_10}
\int\limits_0^1 {f\left( t \right)dt}  = 2 \mathop {\lim }\limits_{L \to \infty } \sum\limits_{\ell  = 1}^L {\sum\limits_{m = 1}^{\left\lfloor {{\frac{M}{2}}} \right\rfloor  + 1} {\frac{1}{{{{\left( {2L} \right)}^{2m - 1}}\left( {2m - 1} \right)!}}} } {\left. {{f^{\left( {2m - 2} \right)}}\left( t \right)} \right|_{t = \frac{{\ell  - 1/2}}{L}}},
\end{equation}
where the notation $\lfloor M/2 \rfloor$ implies the floor function over the ratio $M/2$. Thus, substituting the integrand of the integral \eqref{eq_1} into equations \eqref{eq_10} we obtain
\footnotesize
\begin{equation}\label{eq_11}
\arctan \left( x \right) = 2 \mathop {\lim }\limits_{L \to \infty } \sum\limits_{\ell  = 1}^L {\sum\limits_{m = 1}^{\left\lfloor {{\frac{M}{2}}} \right\rfloor  + 1} {\frac{1}{{{{\left( {2L} \right)}^{2m - 1}}\left( {2m - 1} \right)!}}} } {\left. {\frac{{{\partial ^{2m - 2}}}}{{\partial {t^{2m - 2}}}}\left( {\frac{x}{{1 + {x^2}{t^2}}}} \right)} \right|_{t = \frac{{\ell  - 1/2}}{L}}}.
\end{equation}
\normalsize
Making change of the variable $t \to xt$ in the relation \eqref{eq_6} yields
$$
\frac{{{\partial ^{m - 1}}}}{{\partial {{\left( {xt} \right)}^{m - 1}}}}\left( {\frac{1}{{1 + {x^2}{t^2}}}} \right) = \frac{{{\partial ^m}}}{{\partial {{\left( {xt} \right)}^m}}}\arctan \left( {xt} \right)
$$
or
$$
\frac{{{\partial ^{m - 1}}}}{{\partial {t^{m - 1}}}}\left( {\frac{x}{{1 + {x^2}{t^2}}}} \right) = \frac{{{\partial ^m}}}{{\partial {t^m}}}\arctan \left( {xt} \right)
$$
or
\begin{equation}\label{eq_12}
\frac{{{\partial ^{2m - 2}}}}{{\partial {t^{2m - 2}}}}\left( {\frac{x}{{1 + {x^2}{t^2}}}} \right) = \frac{{{\partial ^{2m - 1}}}}{{\partial {t^{2m - 1}}}}\arctan \left( {xt} \right).
\end{equation}
Consequently, using the relation \eqref{eq_12} we can rearrange the limit \eqref{eq_11} as given by
\footnotesize
\begin{equation}\label{eq_13}
\arctan \left( x \right) = 2 \mathop {\lim }\limits_{L \to \infty } \sum\limits_{\ell  = 1}^L {\sum\limits_{m = 1}^{\left\lfloor {{\frac{M}{2}}} \right\rfloor  + 1} {\frac{1}{{{{\left( {2L} \right)}^{2m - 1}}\left( {2m - 1} \right)!}}{{\left. {\frac{{{\partial ^{2m - 1}}}}{{\partial {t^{2m - 1}}}}\arctan \left( {xt} \right)} \right|}_{t = \frac{{\ell  - 1/2}}{L}}}} }.
\end{equation}
\normalsize

Lastly, according to identity \eqref{eq_8} we have
\footnotesize
\[
\frac{{{\partial ^{2m - 1}}}}{{\partial {t^{2m - 1}}}}\arctan \left( {xt} \right) = \frac{{{{\left( { - 1} \right)}^{2m - 1}}\left( {2m - 2} \right)!{x^{2m - 1}}}}{{2i}}\left( {\frac{1}{{{{\left( {xt + i} \right)}^{2m - 1}}}} - \frac{1}{{{{\left( {xt - i} \right)}^{2m - 1}}}}} \right), \,\, m \ge 1,
\]
\normalsize
and, therefore, from the limit \eqref{eq_13} we obtain
\footnotesize
\[
\begin{aligned}
&\arctan \left( x \right) =\\
& \quad i \mathop {\lim }\limits_{L \to \infty } \sum\limits_{\ell  = 1}^L \sum\limits_{m = 1}^{\left\lfloor {{\frac{M}{2}}} \right\rfloor  + 1} \frac{{{x^{2m - 1}}}}{{2m - 1}} \left( \frac{1}{{{{\left( {x\left( {2\ell  - 1} \right) + 2iL} \right)}^{2m - 1}}}} \right. \left. - \frac{1}{{{{\left( {x\left( {2\ell  - 1} \right) - 2iL} \right)}^{2m - 1}}}} \right)
\end{aligned}
\]
\normalsize
or
\footnotesize
\begin{equation}\label{eq_14}
\begin{aligned}
&\arctan \left( x \right) =\\
& \quad i \mathop {\lim }\limits_{L \to \infty } \sum\limits_{\ell  = 1}^L {\sum\limits_{m = 1}^{\left\lfloor {{\frac{M}{2}}} \right\rfloor  + 1} {\frac{1}{{2m - 1}}\left( {\frac{1}{{{{\left( {\left( {2\ell  - 1} \right) + 2iL/x} \right)}^{2m - 1}}}} - \frac{1}{{{{\left( {\left( {2\ell  - 1} \right) - 2iL/x} \right)}^{2m - 1}}}}} \right)} }.
\end{aligned}
\end{equation}
\normalsize
The equation \eqref{eq_14} for the arctangent function is absolutely identical to \cite{Abrarov2016}
\begin{equation}\label{eq_15}
\arctan \left( x \right) = \mathop {\lim }\limits_{L \to \infty } \sum\limits_{\ell  = 1}^L {\sum\limits_{m = 0}^M {\frac{{{{\left( { - 1} \right)}^m} + 1}}{{{{\left( {2L} \right)}^{m + 1}}\left( {m + 1} \right)!}}} } {\left. {\frac{{{\partial ^m}}}{{\partial {t^m}}}\left( {\frac{x}{{1 + {x^2}{t^2}}}} \right)} \right|_{t = \frac{{\ell  - 1/2}}{L}}}.
\end{equation}
However, unlike equation \eqref{eq_14} this identity requires symbolic programming in order to obtain derivatives up to the order $M$. As a result, the equation \eqref{eq_14} is significantly faster in computation than equation \eqref{eq_15}.

Using the equations \eqref{eq_14} and \eqref{eq_15} for relation
\begin{equation}\label{eq_16}
\pi  = 4\arctan \left( 1 \right)
\end{equation}
we obtain two absolutely identical formulas for the constant pi
\footnotesize
\begin{equation}\label{eq_17}
\begin{aligned}
\pi = 4i \mathop {\lim }\limits_{L \to \infty } \sum\limits_{\ell  = 1}^L {\sum\limits_{m = 1}^{\left\lfloor {{\frac{M}{2}}} \right\rfloor  + 1} {\frac{1}{{2m - 1}}\left( {\frac{1}{{{{\left( {\left( {2\ell  - 1} \right) + 2iL} \right)}^{2m - 1}}}} - \frac{1}{{{{\left( {\left( {2\ell  - 1} \right) - 2iL} \right)}^{2m - 1}}}}} \right)} }
\end{aligned}
\end{equation}
\normalsize
and
\begin{equation}\label{eq_18}
\pi  = 4 \mathop {\lim }\limits_{L \to \infty } \sum\limits_{\ell  = 1}^L {\sum\limits_{m = 0}^M {\frac{{{{\left( { - 1} \right)}^m} + 1}}{{{{\left( {2L} \right)}^{m + 1}}\left( {m + 1} \right)!}}} } {g_{\ell ,m}},
\end{equation}
respectively, where
$$
{g_{\ell ,m}} = {\left. {\frac{{{d^m}}}{{d{t^m}}}\left( {\frac{1}{{1 + {t^2}}}} \right)} \right|_{t = \frac{{\ell  - 1/2}}{L}}}
$$
are the expansion coefficients \cite{Abrarov2016}. The equations \eqref{eq_17} and \eqref{eq_18} provide high-accuracy computation even at relatively small integers $L$ and $M$. In particular, sample computations performed with Wolfram Mathematica (version 9) at enhanced precision mode show that at $L = M = 46$ these equations result in approximated value of pi with $105$ coinciding digits with actual value of pi. However, because of the derivatives up to the order $M$, the limit \eqref{eq_18} requires symbolic programming. As a result, the equation \eqref{eq_17} can provide a significantly faster computation than the equation \eqref{eq_18}.

The relation \eqref{eq_16} is just a simplest example for computation of the constant pi. Since, the integrand of the integral \eqref{eq_1} becomes smoother with decreasing $x$, this facilitates the convergence due to derivatives up to the order $M$. Therefore, we can expect a significant improvement of accuracy in truncation of the equation \eqref{eq_14} as the argument $x$ in the arctangent function $\arctan \left( x \right)$ decreases. As a consequence, the accuracy in computation of pi by equation \eqref{eq_14} can be considerably improved by many orders of the magnitude if instead of equation \eqref{eq_16} we apply, for example, the Gauss formula (see equation (43) in \cite{Lehmer1938})
\begin{equation}\label{eq_19}
\pi  = 4\sum\limits_{n = 1}^9 {{\alpha _n}\arctan \left( {\frac{1}{{{\beta _n}}}} \right)}, \qquad {\frac{1}{\beta_n}} << 1,
\end{equation}
where $\alpha _1  = 2852$, $\alpha _2 =  - 398$,  $\alpha _3 = 1950$, $\alpha _4 = 1850$, $\alpha _5 = 2021$, $\alpha _6 = 2097$, $\alpha _7 = 1484$, $\alpha _8 = 1389$, $\alpha _9 = 808$, $\beta _1 = 5257$, $\beta _2 = 9466$, $\beta _3 = 12943$, $\beta _4 = 34208$, $\beta _5 = 44179$, $\beta _6 = 85353$, $\beta _7 = 114669$, $\beta _8 = 330182$ and $\beta _9 = 485298$. Specifically, applying equation \eqref{eq_14} to the Gauss formula \eqref{eq_19} we have found that at same integers $L = M = 46$ the number of coinciding digits with actual value of the constant pi increases from $105$ to $274$.

Although application of the equation \eqref{eq_14} for computing pi as a sum of the arctangent functions involving only small arguments may be promising, this approach requires further analysis and experimental results that are beyond the scope of the present work.

It should also be noted that due to double summation with respect to the indices $\ell$ and $m$ the application of the equation \eqref{eq_14} may be advantageous for flexible computational optimization by choosing appropriate values of the integers $L$ and $M$.

\section{Conclusion}

A simple identity for the derivatives of the arctangent function is presented as an alternative to the Adegoke--Layeni--Lampret formula. The algorithmic implementation of the proposed identity can significantly accelerate the computation due to no requirement for symbolic programming in determination of the derivatives for the arctangent function.

\section*{Acknowledgments}

This work is supported by National Research Council Canada, Thoth Technology Inc. and York University.



\end{document}